\documentclass[a4paper,12pt]{article}
\usepackage{amsmath}
\usepackage{amsthm}
\usepackage{amssymb}
\usepackage{amscd}
 \makeatletter

 \@addtoreset{equation}{section}
  \makeatother

\theoremstyle{plain}
\newtheorem{thm}{Theorem}[section]
\newtheorem{prop}[thm]{Proposition}
\newtheorem{cor}[thm]{Corollary}
\newtheorem{lem}[thm]{Lemma}
\theoremstyle{definition}
\newtheorem{defn}{Definition}[section]

\theoremstyle{remark}
\newtheorem{rem}{Remark}[section]

\begin{document}
\title{Pre-symplectic structures  on  the space of connections}
\author{ Tosiaki Kori\\
Department of Mathematics\\
School of Science and Engineering\\
 Waseda University \\
3-4-1 Okubo, Shinjuku-ku
Tokyo, Japan.\\ e-mail: kori@waseda.jp
}
\date{ }

\maketitle

\begin{abstract} 
Let \(X\) be a four-manifold with boundary three-manifold \(M\).   We shall describe (i) a pre-symplectic structure on the space \({\cal A}(X)\) of connections on the bundle \(X\times SU(n)\) that comes from the canonical symplectic structure on the cotangent space \(T^{\ast}{\cal A}(X)\), and (ii) a pre-symplectic structure on the space \({\cal A}^{\flat}_0(M)\)  of flat connections on \(M\times SU(n)\)  that have null charge.    These two structures are related by the boundary restriction map.   We discuss also the Hamiltonian features of the space of connections \({\cal A}(X)\) together with the action of the group of gauge transformations.
\end{abstract}


Subj. Class: Global analysis, Quantum field theory.

{\bf Keywords}   Pre-symplectic structures.   Moduli space of flat connections.  Chern-Simons functionals.

\section{Introduction}

Let \(X\) be an oriented Riemannian four-manifold with boundary three-manifold \(M=\partial X\).         For the trivial principal bundle \(P=X\times G\), \(G=SU(n)\), we denote the space of irreducible \(L^2_{2}\)-connections by \(\mathcal {A}(X)\),  \(L^2_l\) being the Sobolev space of square-integrable functions whose distributional derivatives up to order \(l\) are square-integrable.   The tangent space of \({\cal A}(X)\)  is the space of  \(L^2_2\)-valued 1-forms on \(X\):
  \[T_A{\cal A}(X)=\,\mathcal{E}^1_{2}(X,\,ad\,P).\]
  \(ad\,P\) is the bundle of Lie algebra associated to the adjoint representation of \(Lie\,G\).    In this paper we consider only the trivial bundle \(P=X\times SU(n)\), hence 
  \[T_A{\cal A}(X)=\,\mathcal{E}^1_{2}(X,\,Lie\,G).\]
The cotangent space is  the space of  \(L^2_2\)-valued 3-forms on \(X\):
 \[T_A^{\ast}{\cal A}(X)=\,\mathcal{E}^3_{2}(X,\,Lie\,G).\] 
 Correspondingly,  \({\cal A}(M)\) is the space of irreducible \(L^2_{\frac32}\)-connections on \(M\).    By virtue of the Sobolev lemma the boundary restriction map \(r:\mathcal{A}(X)\longrightarrow\mathcal{A}(M)\) is a continuous linear operator.
A connection is said to be flat if its curvature \(F_A=dA+\,A\,\wedge\,A\) vanishes.   
The space of flat connections over \(X\) is denoted by 
\({\cal A}^{\flat}(X)\).    That over \(M\) is denoted by 
\({\cal A}^{\flat}(M)\).    
 The tangent space of \({\cal A}^{\flat}(M)\) is given by
\begin{equation*}
T_A{\cal A}^{\flat}(M)=\{a\in \mathcal{E}^1_{\frac32}(M, Lie\,G ); d_Aa=0\}.
\end{equation*}
We consider the following \(\mathbf{Z}\)-valued functional on \({\cal A}^{\flat}(M)\):
\begin{equation*}
\deg(A)=\frac{1}{24\pi^3}\,\int_M\,Tr\,A^3,\,\qquad A\in{\cal A}^{\flat}(M)\,.
\end{equation*}
 \({\cal A}^{\flat}(M)\) is decomposed into the disjoint union  of connected components:
   \[{\cal A}^{\flat}(M)=\bigoplus_{k\in\mathbf{Z}}\, {\cal A}^{\flat}_k(M)\]
   with 
\begin{equation*}
{\cal A}^{\flat}_k(M)=\left\{A\in{\cal A}^{\flat}(M);\quad \deg(A)=k\,\right\}.
\end{equation*}

For a Riemann surface \(\Sigma\) a symplectic structure on the space of connections over \(\Sigma\times SU(n)\) was introduced by Atiyah-Bott, \cite{AB}, and the geometric quantization theory of the moduli space of flat connections was investigated.    A similar attempt over a four-manifold was carried out in our previous work \cite{K}.
In this paper we shall prove the following theorems.

\begin{thm}\label{presymp1}
Let \(P=X\times SU(n)\) be the trivial \(SU(n)-\)principal bundle on a four-manifold \(X\).   There exists a pre-symplectic structure on the space of irreducible connections \({\cal A}(X)\) that is given by 
the 2-form 
\begin{equation}\label{presympform}
\Omega_A(a,b)=\,\frac{1}{8\pi^3}\int_XTr[(ab-ba)F_A ] -\frac{1}{24\pi^3}\int_{ M}Tr[(ab-ba)A]\,,
\end{equation}
for \(a,b\in T_A{\cal A}(X)\)
\end{thm}

\begin{thm}\label{presymp2}
Let \(\omega\) be a 2-form on \({\cal A}(M)\) defined by 
\begin{equation}
\omega_A(a,b)=-\frac{1}{24\pi^3}\int_{ M}\,Tr[(ab-ba)A]\,,\label{presympform2}
\end{equation}
for \(a,b\in T_A{\cal A}(M)\).    Then \(\left({\cal A}_0^{\flat}(M),\,\omega\vert_{ {\cal A}_0^{\flat}(M)}\,\right)\) is a pre-symplectic manifold.
\end{thm}

The author introduced in \cite{K} the formula in the right hand side of (\ref{presympform}),  so the new feature of Theorem \ref{presymp1} is that this formula is obtained in a canonical manner.   
We explain why we call the 2-form \(\Omega\) {\it canonical}.   On the cotangent bundle \(T^{\ast}{\cal A}(X)\) there exist the canonical 1-form \(\theta\)  and the canonical 2-form \(\sigma=\widetilde d\theta\), where \(\widetilde d\) is the exterior derivative on \({\cal A}(X)\).     Let \(\phi\) be a 1-form on \({\cal A}(X)\).   Then \(\phi\) gives a tautological section of the cotangent bundle \(T^{\ast}{\cal A}(X)\) so that the pullback \(\phi^{\ast}\theta\) of \(\theta\) is characterized by  \(\phi^{\ast}\theta=\phi\).   
The pullback  \(\phi^{\ast}\sigma\) of the canonical 2-form \(\sigma\)  is a closed 2-form on \({\cal A}(X)\).    In particular if we take the 1-form given by the Chern-Simons function \cite{D2}:
\begin{equation*}
\vartheta(A)=\,\frac{1}{24\pi^3} (\,AF_A+F_AA-\frac12 A^3)\in \mathcal{E}^3_{2}(X,\,Lie\,G).
\end{equation*} 
then \(\Omega={\vartheta}^{\ast}\sigma\) is given by equation (\ref{presympform}).      Thus, for a four-manifold \(X\) there is a pre-symplectic form \(\Omega\) on \({\cal A}(X)\) that is induced from the canonical symplectic form \(\sigma\) on the cotangent bundle \(T^{\ast}{\cal A}(X)\) by the generating function \(\vartheta :{\cal A}(X) \longrightarrow T^{\ast}{\cal A}(X)\), which is given by the Chern-Simons function.      We note that, although the Chern-Simons function is generally defined for connections over a three-manifold, we use the same formula for connections \(A\in \mathcal{A}(X)\).   

   Every principal \(G\)-bundle over a three-manifold \(M\) can be extended to a 
principal \(G\)-bundle over a four-manifold \(X\) that cobords \(M\), and for a connection \(A\in {\cal A}(M)\) there is a connection \(\mathbf{A}\in {\cal A}(X)\) that extends \(A\).    So we can detect  a pre-symplectic structure on  \({\cal A}(M)\) as the boundary restriction of the pre-symplectic structure \(\Omega\) on \({\cal A}(X)\).    This does not work in general, that is, the 2-form \(\omega\) of (\ref{presympform2}) on \({\cal A}(M)\) gives only a presymplectic structure twisted by the Cartan 3-form  \(\kappa_A(a,b,c)=-\frac{1}{24\pi^3}\int_{ M}\,Tr[(ab-ba)c]\,\).   But restricted to the space \({\cal A}^{\flat}_0(M)\) of flat connections with degree \(0\), \(\omega\vert_{ {\cal A}_0^{\flat}(M)}\) gives the corresponding pre-symplectic structure.   

Let \({\cal G}_0(X)\) be the group of gauge transformations on \(X\) that are the identity on the boundary \(M\).   Then the closed 2-form \(\Omega\) in (\ref{presympform}) is \({\cal G}_0(X)\)-invariant and the action of   \({\cal G}_0(X)\) on \({\cal A}(X)\) becomes a Hamiltonian action with moment map given by the square of the curvature form  \(F_A^2\), \cite{K}.      Although we know the generalization of the symplectic reduction theory to a pre-symplectic manifold \cite{MR}, we cannot apply this reduction theorem to our case since \(0\) is not a regular value of the moment map \(F_A^2\).   But if the boundary \(M\) is not empty, the space of flat connections \({\cal A}^{\flat}(X)\) is a smooth manifold contained in the zero level set.   Hence the pre-symplectic structure (\ref{presympform}) on \({\cal A}(X)\) descends to the moduli space of flat connections \({\cal M}^{\flat}(X)={\cal A}^{\flat}(X)/{\cal G}_0(X)\).   The pre-symplectic structure on \({\cal M}^{\flat}(X)\) is given by the restriction of \(\Omega\) to the space of flat connections:
\begin{equation*}
\omega^{\flat}_{[A]}([a]\,,\,[b]\,)=-\frac{1}{24\pi^3}\int_{ M}\,Tr[(ab-ba)A]\,
\end{equation*}
for 
\([A]\in {\cal M}^{\flat}(X)\) and \([a],\,[b]\in T_{[A]}{\cal M}^{\flat}(X)\), with 
\(a,b\in T_A{\cal A}^{\flat}(X)\).

On the other hand we shall prove in section 4 the following  
\begin{thm}    
The boundary restriction map
\begin{equation*}
 \overline r: {\cal M}^{\flat}(X)\longrightarrow  {\cal A}^{\flat}_0(M)
\end{equation*}
is a diffeomorphism.
\end{thm}
Hence we have an isomorphism of pre-symplectic manifolds
\begin{equation*}
\left({\cal M}^{\flat}(X),\,\omega^{\flat}\,\right)\stackrel{\simeq}{\longrightarrow}\,\left({\cal A}^{\flat}_0(M),\,\omega\,\right).
\end{equation*}

\section{ Preliminaries on the space of connections}

\subsection{Calculations on the space of connections}

Let \(M\) be a compact, connected and oriented \(m\)-dimensional Riemannian manifold  possibly with boundary \(\partial M\).    Let \(P\stackrel{\pi}{\longrightarrow}M\) be a principal \(G\)-bundle, \(G=SU(N)\), \(N\geq 2\).   A connection 1-form on \(P\) is a \(Lie\,G\)-valued 1-form on \(P\) that is invariant under \(G\), acting by a combination of the action on \(P\) and the adjoint action on \(Lie\,G\).     
 The connections over \(P\) form an affine space modeled on  the vector space of \(Lie \,G\)-valued 1-forms on \(M\), and the group of automorphisms on \(P\) acts on it.   
We shall deal with the quotient space of the space of connections by the action of the automorphism group.       In order to show that there is a good quotient space,  we need to show that the action has local slices.   To establish such a result we adopt the analysis in terms of Sobolev completions.     References  for this are 
 \cite{ DK, FU, KN, M1, M2} etc.

 We introduce the \(L^2_l\)-norms on the space of smooth sections of a vector bundle \(V\longrightarrow M\) with a positive definite inner product:
 \[\Vert s\Vert_{L^2_l}\,=\,\int_M\,\left(\sum_{i=0}^l\,\Vert \nabla^i(s(x)\Vert^2\right)^{\frac12}\,d{\rm vol}\,.\]
Here \(\nabla\) is the covariant derivative with respect to a connection on \(V\) and the Levi-Civita connection on \(T^{\ast}M\).    The norms in the integrand are the point-wise norms induced by the inner products on \(T^{\ast}M\) and \(V\).
 For each non-negative integer \(l\),  \(L^2_l(M, V)\) is the completion of the space of smooth sections with respect to this norm, and is called the Sobolev space of \(L^2_l\)-sections.    For \(V=\wedge^r\,T^{\ast}M\) we have  the Sobolev space \(\mathcal{E}^r_l(M)\) of r-forms on \(M\).    In the sequel we assume 
 \[l\geq1+\frac{m}{2}\,.\]
    The Sobolev multiplication theorem implies that there is a continuous linear multiplication mapping 
 \[\mathcal{E}^p_l(M)\,\otimes\,\mathcal{E}^q_l(M)\,\longrightarrow \mathcal{E}^{p+q}_l(M)\,.\]
 Moreover the definition can be extended to the \(Lie\,G\) valued differential forms on \(M\).   Let \(\mathcal{E}^r_l(M,\,Lie\,G)\) be the Sobolev space of \(Lie\,G\)-valued r-forms on \(M\).    It is the completion of \(\{\phi\otimes\,X\,;\,\phi\in \mathcal{E}^r_l(M),\,X\in \,Lie\,G\,\}\). 
 
 We define the dual pairing of  \(\phi\in \mathcal{E}^{p}_l(M)\, \) and \(\psi\in\mathcal{E}^{m-p}_l(M),\)  by
\begin{equation}\label{dualpair}
\langle\,\phi\,,\,\psi\,\rangle\,=\,\int_M\,\phi\wedge \,\psi\,.
\end{equation}
It is a continuous bilinear form on \( \mathcal{E}^{p}_l(M)\, \times \mathcal{E}^{m-p}_l(M)\).
   We have also the dual pairing of 
\(\mathcal{E}^p_l(M,\,Lie\,G)\) and \(\mathcal{E}^{m-p}_l(M,\,Lie\,G)\), where the dual of \(Lie\,G\) is identified with \(Lie\,G\) by the non-degenerate symmetric bilinear form \((X,Y)\longrightarrow tr(\,XY\,)\).     (\ref{dualpair}) is extended to 
the \(Lie\,G\)-valued forms:
\begin{equation}\label{dualpair2}
\langle\,\phi\otimes X\,,\,\psi\otimes Y\,\rangle\,=\,\langle\,\phi\,, \,\psi\,\rangle\, tr\,(XY)\,, 
\end{equation}
for \(\phi\in \mathcal{E}^{p}_l(M)\, \), \(\psi\in\mathcal{E}^{m-p}_l(M),\) and \(X,\,Y\in Lie\,G\).

 Let \(Ad\,P\) be the fiber bundle \(Ad\,P=P\times_GG\longrightarrow\,M\).   We know that the group of smooth sections of \(Ad\,P\) under fiber-wise multiplication is isomorphic to \(Aut\,(P)\);   \(Aut\,P\) is the group of \(G\)-equivariant diffeomorphisms of \(P\) that cover the identity on the base \(M\).    It is called  the {\it group of gauge transformations}.   This group is an infinite dimensional Lie group and its Lie algebra 
is the associated vector bundle  \(\,ad\,P\,=\,P\times_G\,Lie\,G\) obtained by the adjoint representation of \(G\) on \(Lie\,G\).

  We write \({\cal A}={\cal A}(M)\) the space of {\it irreducible} connections over \(P\) modeled on the vector space \(\mathcal{E}^1_{l-1}(M,\,ad\,P)\).   
So the tangent space at \(A\in{\cal A}\) is
 \begin{equation*}
 T_A{\cal A}=\mathcal{E}^1_{l-1}(M, ad \,P)\,.
 \end{equation*}
 A connection \(A\) is  irreducible if the holonomy group of \(A\) is an open subgroup of \(G\).    The irreducibility assumption is necessary to fix the  horizontal distribution of \(T\mathcal{A}\) associated to the connection ( Coulomb gauge condition \cite{S}).   In setions 1.2 and 1.3 we shall give a short explanation.      
 
The cotangent space at \(A\in\mathcal{A}\) is
  \begin{equation*}
 T_A^{\ast}{\cal A}=\mathcal{E}^{m-1}_{l-1}(M, ad\,P)\,).
 \end{equation*}
The dual pairing of \( T_A{\cal A}=\,\mathcal{E}^1_{l-1}(M, ad\,P)\) and \(T^{\ast}_A{\cal A}=\mathcal{E}^{m-1}_{l-1}(M, ad\,P)\,)\) is given by (\ref{dualpair2}) :
\[\langle \alpha,\,a\rangle_A=\int_M\,tr(\,\alpha\wedge a\,)\,.\quad a\in T_A{\cal A}\,,\,\alpha\in T^{\ast}_A{\cal A}\,. \] 

We write the group of \(L^2_l\)-gauge transformations by \({\cal G}^{\prime}(M)\):
\begin{equation*}
{\cal G}^{\prime}(M)=\,\mathcal{E}^0_l(M, Ad\,P)\, .\end{equation*}
  \({\cal G}^{\prime}(M)\) acts on \({\cal A}(M)\) by 
\begin{equation*}
g\cdot A=g^{-1}dg+g^{-1}Ag=A+g^{-1}d_Ag.
\end{equation*}
By the Sobolev lemma one sees that \({\cal G}^{\prime}(M)\) is a Hilbert Lie\,group  and its action is a smooth map of Hilbert manifolds.   

\begin{rem}
Instead of \(\mathcal{A}(M)\) and  \(\mathcal{G}(M)\), several authors write \(\mathcal{A}(P)\) and \(\mathcal{G}(P)\), for example \cite{M2}.    We prefer our notation because later we shall deal with connections over the four-manifold \(X\) as well as over the boundary three-manifold \(M=\partial\,X\) at the same time so that we have  \(\mathcal{A}(X)\) and \(\mathcal{A}(M)\), respectively \(\mathcal{G}(X)\) and \(\mathcal{G}(M)\).   We prefer to avoid notation such as \(\mathcal{A}(P\vert M)\), etc.
\end{rem}

In the following we choose a fixed point \(p_0\in M\) and deal with the group of gauge transformations that are the identity at \(p_0\):
\[{\cal G}={\cal G}(M)
=\{g\in {\cal G}^{\prime}(M);\,g(p_0)=1\,\}.\] 
  \({\cal G}\) acts freely on \({\cal A}\).    
 Let \({\cal C}(M)={\cal A}(M)/{\cal G}(M)\) be the quotient space of this action.   It is a smooth  manifold.     
 
   We have 
    \[Lie\,({\cal G})=\mathcal{E}^0_l(M, ad\,P)\,.\]

The action of \(Lie\,({\cal G})\) on \({\cal A}\) is given by the covariant exterior derivative:
\begin{equation*}
d_A=d+[A\,,\quad]:\, Lie\,\mathcal{G}=\mathcal{E}^0_l(M, ad\,P)\longrightarrow \mathcal{E}^1_{l-1}(M, ad\,P)=T_A\mathcal{A}\,.
\end{equation*}

The fundamental vector field on \({\cal A}\) corresponding to  
\(\xi \in Lie({\cal G})\) is given by 
\begin{equation*}
d_A\xi\,=\,\frac{d}{dt}\lvert_{t=0}(\exp \,t\xi)\cdot A\,.
\end{equation*}
Since \(A\in \mathcal{A}\) is irreducible we have  \(Ker\,d_A\,=\,0\)  .     
The tangent space to the orbit at \(A\in{\cal A}\) is 
\begin{equation*}
T_A({\cal G}\cdot A)=\{d_A\xi\,;\,\xi\in \mathcal{E}^0_l(M, ad\,P)\}.
\end{equation*}
\\

Now we explain some calculations on the affine space \(\mathcal{A}\).    We shall look at the formulae for the differentiation of functions, vector fields and differential forms \cite{BN, BV, DK, MW, MV}.   
 A vector field  \(\mathbf{a}\) on \({\cal A}\) is a section of the tangent bundle;
 \(\mathbf{a}(A)\in T_A{\cal A}\) for \(\forall A\in \mathcal{A}\).   A 1-form \(\varphi\) on \({\cal A}\) is a section of the cotangent bundle;
 \(\varphi(A)\in T_A^{\ast}{\cal A}\) for \(\forall A\in \mathcal{A}\).       Similarly we consider any \(r\)-form \(\psi\); \(\psi(A)\in \wedge^rT^{\ast}_A\mathcal{A}\,\).
    
For a smooth function \(F=F(A)\) on \({\cal A}\) taking the value in a vector space \(V\) the derivative \(\partial F(A)\) is defined by the functional derivation with respect to \(A\in {\cal A}\):
\begin{equation*}
(\partial F(A))a\,=\,\lim_{t\longrightarrow 0}\frac{1}{t}\left(\,F(A+ta)-F(A)\,\right),\quad \mbox{ for \(a\in T_A{\cal A}\)}.
\end{equation*}
 Hence \(\partial F(A)\in T^{\ast}_A\mathcal{A}\otimes V\).
 
For example, let  \(f( A)= A-A_0 \in \mathcal{E}^1_{l-1}(M,ad\,P)\)
 for a fixed connection \(A_0\in \mathcal{A}\).   Then 
\[(\partial f (A))\,a=a.\]

The curvature of a connection \(A\in {\cal A}\) is by definition 
\[F(A)=dA+\frac12[\,A\,,\, A\,]\in \mathcal{E}^2_{l-2}(M, ad\,P).\]  
Thus
\[F(A+a)=F(A)+d_Aa+a\wedge a,\]
and we have 
  \[(\partial F(A))a=d_Aa.\]
  
  The derivative of a vector field \(\mathbf{v}\) on \({\cal A}\) and that of a 1-form \(\varphi\) are defined similarly:
  \[(\partial \mathbf{v}(A)) \in T^{\ast}_A\mathcal{A}\otimes T_A{\cal A}\,, \qquad (\partial \varphi(A))\in T^{\ast}_A{\cal A}\otimes T^{\ast}_A{\cal A}\,, \]
 that is,
  \[(\partial \mathbf{v}(A))a \in T_A{\cal A},\qquad(\partial \varphi(A))a \in T^{\ast}_A{\cal A}\,,\quad\mbox{ for } \forall a\in T_A{\cal A}.\]
   It follows that the derivative of a function   \(F=F(A)\) with respect to a vector field \(\mathbf{v}\) is given by
  \[(\mathbf{v}F)(A)=(\partial F(A))\,v ,\] 
  where \(v=\mathbf{v}(A)\in T_A\mathcal{A}\).   
  
    We have the following formulas.
  \begin{equation}
 [\, \mathbf{v},\,\mathbf{w}\,](A)=(\partial \mathbf{v}(A))w-(\partial\mathbf{w}(A))v\,,\label{derivation1}
 \end{equation}
 \begin{equation}
(\partial<\varphi, \mathbf{u}>(A))\,v\,=\,
 \langle\, \varphi(A),\, (\partial\mathbf{u}(A))v\,\rangle +\langle \,(\partial\varphi(A))\,v\,,\,u\,\rangle ,\label{derivation2}
  \end{equation}
    where \(u=\mathbf{u}(A)\in T_A\mathcal{A}\), etc.

Let \(\widetilde d\)  be the exterior derivative on \({\cal A}\).   For a function \(F\) on \({\cal A}\), the exterior derivative  \((\widetilde d\,F)_A\in T^{\ast}_A\mathcal{A}\) of \(F\) at the point \(A\in\mathcal{A}\) is defined by \((\widetilde d\,F)_A\,a=(\partial F(A))\,a\) for any  \( a\in T_A\mathcal{A}\).   \\
 For a 1-form \( \varphi\) on \({\cal A}\), 
\begin{eqnarray}
(\widetilde d \varphi)_A( a\,,\, b\,)&=&(\partial< \varphi,{\bf b}>(A)) a\,-\,(\partial< \varphi,{\bf a}>(A)) b\,-< \varphi,\, [{\bf a},{\bf b}]>(A)\nonumber \\[0.2cm]
&=& <(\partial \varphi(A))\, a\,,  \, b\,>-<(\partial \varphi(A))\, b\,,\, a\,>,\label{extder1}
\end{eqnarray}
where \(a,\,b\in T_A\mathcal{A}\) and \(\,\mathbf{a}\,,\,\mathbf{b}\) are vector fields such that \(\mathbf{a}(A)=a\) and  \(\,\mathbf{b}(A)=b\).   In fact,
the first line is the definition of the exterior derivative \(\widetilde d\, \varphi\), \cite{KN}, and (\ref{extder1}) follows from (\ref{derivation1}) and (\ref{derivation2}).
 Likewise, if \(\varphi\) is a 2-form on \({\cal A}\), 
\begin{eqnarray*}
(\widetilde d\,\varphi)_A( a,\, b\,,\, c\,)&=&\,<(\,\partial\varphi(A))( b, c),\,a\,>\,+\,<(\partial\varphi(A))(c,a)\,,\,b\,>\\
\,&\,&\,+\,<\,(\partial\varphi (A))( a,\, b)\,,\,c\,>.
\end{eqnarray*}

\subsection{The moduli spaces \({\cal A}/{\cal G}\) and \({\cal A}/{\cal G}_0\)}

Let \(\mathcal{G}^{\prime}(M)\) be ( as in the preceding paragraph ) the group of \(L^2_l\)-gauge transformations on the trivial principal \(G\)-bundle \(P=M\times G\).    When \(M\) has boundary we must deal with the boundary behavior of gauge transformations.     Here we shall  mention two types: the Dirichlet boundary condition and the Neumann boundary condition.   Let \({\cal G}^{\prime}(\partial M)\) be 
the group of \(L^2_{l-\frac12}\) gauge transformations on the boundary \(\partial M\,\).   We have the restriction map to the boundary:       
\begin{equation*} 
r:\, {\cal G}^{\prime}(M) =\mathcal{E}^0_l(M,G)\,
\longrightarrow \, {\cal G}^{\prime}(\partial M)=\,\mathcal{E}^0_{l-\frac{1}{2}}(\partial M,G)\,.
\end{equation*}
\(r\) is  a continuous map.   Moreover it is a compact map by the Sobolev embedding theorem.
Let \({\cal G}_0={\cal G}_0(M)\) be the kernel of the restriction map.   
It is the group of gauge transformations that are the identity on the boundary.
\({\cal G}_0\) acts freely on \({\cal A}=\mathcal{A}(M)\) and \({\cal A}/{\cal G}_0\) is therefore a smooth infinite dimensional manifold.   

In the following we choose a fixed point \(p_0\in M\) on the boundary \(\partial M\) and deal with the group of gauge transformations based at \(p_0\):
\[{\cal G}={\cal G}(M)
=\{g\in {\cal G}^{\prime}(M);\,g(p_0)=1\,\}.\] 
If \(\partial M=\phi\), \(p_0\) may be any point of \(M\).   \({\cal G}\) acts freely on \({\cal A}\) and the orbit space  \({\cal A}/{\cal G}\) is a smooth infinite dimensional manifold.    \\
Let
 \[{\cal G}(\partial M)=\{g\in {\cal G}^{\prime}(\partial M);\,g(p_0)=1\,\}\,.\]
    \({\cal G}_0\) is the kernel of the restriction map \(r:\,{\cal G}(M)\longrightarrow {\cal G}(\partial M)\).   
  We have
\begin{equation*}
Lie ({\cal G}_0)=\{\xi\in Lie({\cal G});\,\xi\vert\partial M=0\}.
\end{equation*}

The restriction to the boundary yields also the following map on the connection spaces:
\begin{equation*} 
r:\, {\cal A}(M) = \mathcal{E}^1_{l-1}(M,G)\,
\longrightarrow \, {\cal A}(\partial M)=\,\mathcal{E}^1_{l-\frac{3}{2}}(\partial M,G)\,.
\end{equation*}

We introduce two moduli spaces of irreducible connections;
\begin{equation*}
 {\cal B}\,=\, {\cal B}(M)={\cal A}/{\cal G}_0 ,\qquad {\cal C}\,=\,{\cal C}(M)={\cal A}/{\cal G}.
\end{equation*}

\({\cal B}\) is a \({\cal G}/{\cal G}_0\)-principal bundle over \({\cal C}\).    \({\cal C}\) coincides with \({\cal B}\) if  \(M\) has no boundary.     \({\cal C}\) is finite dimensional but in general \({\cal B}\) is infinite dimensional, in fact it contains the orbit of \({\cal G}(\partial M)\).

\(  {\cal B}(M)\) is a smooth manifold modeled locally on the ball in the subspace \(\ker d^{\ast}_A\) of the Hilbert space \(\mathcal{E}^1_{l-1}(M, Lie\,G)\).       \({\cal C}(M)\)  is a smooth manifold modeled locally on the ball in the Hilbert subspace \( \ker d^{\ast}_A \cap \ker(\ast\vert\partial M)\subset\,\mathcal{E}^1_{l-1}(M, Lie\,G)\).     The reader can find precise and technical descriptions of these facts in ~\cite{BV, D, DK, M2,  MV}.       In the following we shall give a short explanation about  connections on the principal bundles \({\cal A}\longrightarrow {\cal B}\) and 
\({\cal A}\longrightarrow {\cal C}\).    

The Stokes formula is stated as follows:
\begin{equation}\label{stokes}
\int_{\partial M}\,<\,f,\, \ast u\,>=\int_M<\,d_Af,\, \ast u\,>\,-\, \int_M<\,f,\,\ast d^{\ast}_Au\,>,
\end{equation}
for \(f\in\mathcal{E}^0_l(M,Lie\,G),\, u\in\mathcal{E}^1_{l-1}(M,Lie\,G)\).
If  \(M\) is a compact manifold without boundary  (\ref{stokes}) implies the following decomposition:
\begin{equation*}
T_A{\cal A}=\{d_A\xi;\,\xi\in Lie({\cal G})\,\}\,\oplus \,H^0_A,\end{equation*}
where
\[H^0_A=\{a\in\mathcal{E}^1_{l-1}(M, Lie\,G);\, d_A^{\ast}a=0\}.\]
In this case we have
\begin{equation*}
T_{[A]}{\cal B}=T_{[A]}{\cal C}\simeq H^0_A,
\end{equation*}
where \([A]\in {\cal B}\) is the equivalence class of \(A\in\mathcal{A}\).

We shall look at the case when \(M\) has boundary.
Let \(\Delta_A\) be the covariant Laplacian defined as the closed extension of \(d^{\ast}_Ad_A \) with the domain of definition 
\({\cal D}_{\Delta_A}=\{u\in \mathcal{E}^0_l(M,Lie\,G);\,u\vert\partial M=0\}\) .   
Since \(A\in{\cal A}\) is irreducible \(\Delta_A:{\cal D}_{\Delta_A}\longrightarrow \mathcal{E}_{l-2}^0(M,Lie\,G)\) is an isomorphism.  Let 
\(G_A=(\Delta_A)^{-1}\) be the Green operator of the Dirichlet problem :
\[\left\{\begin{array}{ccc} 
\Delta_A \,u&=& f\\[0.2cm]
u\vert \partial M &=& 0\end{array}\right.\]

Let \(A\in{\cal A}\).   
By (\ref{stokes}) we have the orthogonal decomposition:
\begin{equation*}
T_A{\cal A}=\{d_A\xi;\,\xi\in Lie({\cal G}_0)\,\}\,\oplus \,H^0_A,
\end{equation*}
with
\[H^0_A=\{a\in\mathcal{E}^1_{l-1}(M, Lie\,G);\, d_A^{\ast}a=0\}.\]
\(a\in T_A{\cal A}=\mathcal{E}^1_{l-1}(M, Lie\,G)\) is decomposed to 
\begin{equation}\label{Ddecomp}
a=d_A\xi+b,\quad\mbox{ with }\, \xi=G_Ad^{\ast}_Aa\in Lie({\cal G}_0),\quad b\in H^0_A.
\end{equation}

From this we see that the \({\cal G}_0\)-principal bundle \(\pi:{\cal A}\longrightarrow {\cal B}\) has a natural connection defined by the horizontal subspace \(H^0_A\), 
which is given by the connection 1-form
\(\gamma^0_A=G_Ad_A^{\ast}\).    
The curvature form of 
\(\gamma^0\) is given by
\[{\cal F}^0_A(a,b)=G_A(\ast\,[a,\,\ast b])
\qquad\mbox{ for \(a,b \in H^0_A\),}\]
\cite{D, S}.\\

Now we proceed to the fibration  \({\cal A}\longrightarrow {\cal C}={\cal A}/{\cal G}\).   

Let \(\Delta_A^{(n)}\) be the closed extension of \(d^{\ast}_Ad_A \) with the domain of definition 
\({\cal D}_{\Delta^{(n)}_A}=\{u\in \mathcal{E}^0_l(M,\,Lie\,G)\,;\,\,\ast d_Au\vert\partial M=0\}\) .    ( The superscript \((n)\) indicates the Neuman boundary condition.)

For a 1-form \(v\) on \(M\), let \(g=\,R_Av\in \mathcal{E}^0_{l}(M,\,Lie\,G)\) denote the solution of the following boundary value problem:
\[\left\{
\begin{array}{ccc} \Delta^{(n)}_A \,g &=& 0\qquad \\[0.2cm]
\ast d_A g\vert\partial M&=&\ast v \vert\partial M.
\end{array}\right.\]

Let \(A\in{\cal A}\).     
We have the orthogonal decomposition:
\begin{equation*}
T_A{\cal A}=\{d_A\xi;\,\xi\in Lie({\cal G})\}\,\oplus \,H^{(n)}_A,\end{equation*}
where
\[H^{(n)}_A=\{a\in\mathcal{E}^1_{l-1}(M, Lie\,G);\, \,d_A^{\ast}a=0,\,\mbox{and } \,\ast a\vert\partial M=0\}.\]

In fact let \(a\in\mathcal{E}^1_{l-1}(M,Lie\,G)\) and \(a=d_A\xi+b\) be the decomposition of 
(\ref{Ddecomp}
), then \(\xi=\,G_Ad^{\ast}_Aa\in Lie({\cal G}_0)\,\) and \(b\in H^0_A\).   Put \(\eta=R_Ab\).   
Then we have the orthogonal 
decomposition 
\[a=d_A(\xi+\eta)+c,\]
with \(c\in H^{(n)}_A\) and \(\xi+\eta\in Lie({\cal G})\).   

If we write
\begin{equation*}
\gamma_A=\gamma^0_A+R_A(I-d_A\gamma^0_A) =\, N_Ad_A^{\ast}\,+\,R_A,
\end{equation*}
where \(I\) is the identity transformation on \(T_A{\cal A}\), 
then \(\gamma_A\) is a \(Lie({\cal G})-\)valued 1-form which vanishes on \(H^{(n)}_A\) and \(\gamma_A d_A\xi=\xi\),   
that is, \(\gamma_A\) is the connection 1-form of the fibration  \({\cal A}\longrightarrow {\cal C}\).    
The curvature form of \(\gamma_A\) is given by
\begin{equation*}
{\cal F}_A(a,b)=N_A(\ast[a,\,\ast b])\qquad\mbox{for \(a,b\in H_A\)}.
\end{equation*}
Here  \(N_A=(\Delta_A^{(n)})^{-1}\) is the Green operator of Neuman problem:
\[\left\{
\begin{array}{ccc} \Delta_A^{(n)} \,v &=&f\qquad\\[0.2cm]
\ast d_Av \vert\partial M&=& 0\qquad\mbox{on \(\partial M\)}.
\end{array}\right.\]

\subsection{ The moduli space of flat connections }

 In the sequel we shall suppress the Sobolev indices.     Thus \({\cal A}\) is always the space of irreducible \(L^2_{l-1}\) connections and \({\cal G}\) is the group of based \(L^2_l\) gauge transformations.   

The space of flat connections is 
\begin{equation*}
{\cal A}^{\flat}(M)=\{A\in{\cal A}(M); \, F_A=0\},
\end{equation*}
which we shall often abbreviate to \({\cal A}^{\flat}\).   
The tangent space of \({\cal A}^{\flat}\) is given by
\begin{equation*}
T_A{\cal A}^{\flat}=\{a\in \mathcal{E}^1_{l-1}(M, Lie\,G);\, d_Aa=0\}.
\end{equation*}

The moduli space of flat connections is by definition 
\begin{equation}\label{flatmoduli}
{\cal M}^{\flat}={\cal A}^{\flat}/{\cal G}_0 .\end{equation}
If there is any doubt about which manifold is involved,  we shall write 
\({\cal M}^{\flat}( M)\) for the orbit space \({\cal M}^{\flat}={\cal A}^{\flat}(M)/{\cal G}_0(M)\).   

We know that  \({\cal M}^{\flat}\) is a smooth manifold.    In fact, the coordinate mappings are described by the implicit function theorem~\cite{ DK, FU, M2, MV}.     
  
For \(A\in{\cal A}^{\flat}\) there is a slice for the \({\cal G}_0\)-action on \({\cal A}^{\flat}\) given by the Coulomb gauge condition:
\begin{equation*}
V_A= \{a\in \mathcal{E}^1_{l-1}(M,\,Lie\,G); \,\vert a\vert \,< \epsilon,\,d_Aa+a\wedge a=0,\,d^{\ast}_Aa=0\}.
\end{equation*}  
Let 
\begin{equation}\label{flattangentsp}
 H^{\flat}_A=\{\mathcal{E}^1_{l-1}(M, Lie\,G);\, \,d_Aa=0,\,d^{\ast}_Aa=0\}.
 \end{equation}
The Kuranishi map is defined by 
\[K_A:\, \mathcal{E}^1_{l-1}(M,  Lie\,G)\ni \alpha\longrightarrow K_A(\alpha)=\alpha+d_A^{\ast}G_A(\alpha\wedge\alpha)\in\mathcal{E}^1_{l-1}(M,  Lie\,G).\]
Since the derivative of \(K_A\) at \(\alpha=0\)  becomes the identity transformation on \(\mathcal{E}^1_{l-1}(M,  Lie\,G)\), the implicit function theorem in Banach space yields that  \(K_A\) gives an isomorphism on a small neighborhood of \(0\).   Thus we see that 
the slice \(A+V_A\) is a neighborhood of \(A\) that is homeomorphic to the following subset of \(H^{\flat}_A\);
\[ \{\beta\in H^{\flat}_A;\,\vert \beta\vert \,<\epsilon,\, \lambda_A(\beta)=0\},\]
where
\begin{equation*}
\lambda_A(\beta)=(I-G_A\Delta_A)(\alpha\wedge\alpha),\quad \alpha=K_A^{-1}\beta.   
\end{equation*}
 \hfill\qed 
 \\

We investigate also the moduli space of flat connections modulo the total gauge transformation group \({\cal G}\), 
\begin{equation}\label{flatmodulin}
{\cal N}^{\flat}={\cal A}^{\flat}/{\cal G}.
\end{equation}
A slice in a neighborhood of \(A\in {\cal A}^{\flat}\) in this case is 
\begin{equation*}
W_A= \{a\in \mathcal{E}^1_{l-1}(M,\, Lie\,G); \,\vert a\vert \, < \epsilon,\,d_Aa+a\wedge a=0,\,d^{\ast}_A=0,\mbox{and}\, \ast a\vert{\partial M}=0 \}.
\end{equation*}  
The Kuranishi map is defined by 
\[L_A:\, \mathcal{E}^1_{l-1}(M,  Lie\,G)\ni \alpha\longrightarrow L_A(\alpha)=\alpha+d_A^{\ast}N_A(\alpha\wedge\alpha)\in\mathcal{E}^1_{l-1}(M,  Lie\,G).\]
The same argument as above yields that there is a slice through \(A\) in \({\cal N}^{\flat}\) that is homeomorphic to 
\[ \{\beta\in H^{\flat}_A;\,\vert \beta\vert \,<\epsilon,\, \mu_A(\beta)=0\},\]
where
\[\mu_A(\beta)=(I-N_A\Delta_A)(\alpha\wedge\alpha),\quad \alpha=L_A^{-1}\beta.\]

The dimension of \({\cal N}^{\flat}\) is finite but  the dimension of \({\cal M}^{\flat}\) is in general  infinite.

\section{Canonical structure on \(T^{\ast}{\cal A}\)}

The cotangent bundle of any manifold is endowed with the canonical symplectic structure and we have the standard theory of Hamiltonian mechanics, \cite{BV, M1, MR}.   Here we apply these standard techniques of symplectic geometry to our infinite dimensional manifold \({\cal A}(M)\) and give some explicit formulas.    

\subsection{ Canonical 1-form and 2-form on \(T^{\ast}{\cal A}\)}

Let \(M\) be a manifold of \(\dim M=m\) possibly with the non-empty boundary.   Let \(P=M\times G\) be the trivial \(G\)-bundle.   Let 
 \({\cal A}={\cal A}(M)\) the space of irreducible \(L^2_{l-1}\) connections, \(l >1+\frac m2\).   
Let \(T^{\ast}{\cal A}\,\stackrel{\pi}{\longrightarrow}{\cal A}\,\) be the cotangent bundle.    Any point of \(T^{\ast}{\cal A}\) is denoted by \((A,\lambda)\in T^{\ast}{\cal A}\) with \(A\in \mathcal{A}\) and \(\lambda\in T^{\ast}_A{\cal A}\).   
The tangent space  to  the cotangent space \(T^{\ast}{\cal A}\) at the point \((A,\lambda)\in T^{\ast}{\cal A}\) becomes 
\[T_{(A,\lambda)}T^{\ast}{\cal A}=T_A{\cal A}\oplus T^{\ast}_A{\cal A}=\mathcal{E}^1_{l-1}(M,\,Lie\,G)\oplus\mathcal{E}_{l-1}^{m-1}(M.\,Lie\,G).\]
The canonical 1-form \(\theta\) on the cotangent space is defined by:      
\begin{equation}\label{canonical1}
\theta_{(A,\lambda)}(\left(\begin{array}{c}a\\ \alpha\end{array}\right))=\langle\, \lambda,\pi_{\ast}\left(\begin{array}{c}a\\ \alpha\end{array}\right)\,\rangle_A=\langle\, \lambda,\,a \,\rangle_A\,=\,    \int_M\,tr \,a\wedge \lambda\,,
\end{equation}
for all 
\(\left(\begin{array}{c}a\\ \alpha\end{array}\right)\in T_{(A,\lambda)}T^{\ast}{\cal A}\,\).

\begin{lem}\label{lem1-1}
The derivative of the 1-form \(\theta\,\)
 is given by
\begin{equation*}
(\partial \theta(A,\lambda))\, \left(\begin{array}{c}a\\ \alpha\end{array}\right) = \,\langle \alpha\,,\,a\,\rangle\,,\quad 
\forall \left(\begin{array}{c}a\\ \alpha\end{array}\right)\in T_{(A,\lambda)}T^{\ast}{\cal A}\, .
\end{equation*}
\end{lem}

In fact, 
\[
(\partial \theta(A,\lambda))\, \left(\begin{array}{c}a\\ \alpha\end{array}\right) = \lim_{t\longrightarrow 0}\frac{1}{t}\int_M\,(\,tr\,a\wedge (\lambda+t\alpha)\,-\,tr\,a\wedge \lambda\,)
=
\int_M\,tr\,a\wedge \alpha
.
\]

The canonical 2-form \(\sigma\) on the cotangent space is defined by:
\begin{equation}\label{canonical2}
\sigma=\widetilde{d}\theta
.\end{equation}

 Lemma \ref{lem1-1}  and  (\ref{extder1}) yield the following
\begin{prop}
\begin{eqnarray*}
\sigma_{(A,\lambda)}\,(\left(\begin{array}{c}a\\ \alpha\end{array}\right),\,\left(\begin{array}{c}b\\ \beta\end{array}\right)\,)&=&\int_M\,tr[ \,b\wedge\alpha-a\wedge\beta\,]\\[0.2cm]
\quad 
&\forall &\left(\begin{array}{c}a\\ \alpha\end{array}\right)\,,\,\left(\begin{array}{c}b\\ \beta\end{array}\right)\in T_{(A,\lambda)}T^{\ast}{\cal A}\, .
\end{eqnarray*}
\end{prop}

\(\sigma\) is a non-degenerate closed 2-form on the cotangent space \(T^{\ast}{\cal A}\).     We see the non-degeneracy as follows.   Let \(\left(\begin{array}{c}a\\ \alpha\end{array}\right)\in T_{(A,\lambda)}T^{\ast}{\cal A}\), then \(a\in\mathcal{E}_{l-1}^1(M,\,Lie\,G)\) and \(\alpha\in \mathcal{E}_{l-1}^{m-1}(M,Lie\,G)\).   Hence 
\(\ast\alpha\in\mathcal{E}_{l-1}^1(M.Lie\,G)\) and \(\ast a\in\mathcal{E}_{l-1}^{m-1}(M,\,Lie\,G)\) and we have 
\[\sigma_{(A,\lambda)}\,(\left(\begin{array}{c}a\\ \alpha\end{array}\right),\,\left(\begin{array}{c}\ast \alpha\\ \ast  a\end{array}\right)\,)=||\alpha||^2-||a||^2.\]
The formula implies the non-degeneracy of \(\sigma\).\\

Let \(\Phi=\Phi(A,\lambda) \) be a function on the cotangent space \(T^{\ast}{\cal A}\).    
The corresponding Hamiltonian vector field 
\(X_{\Phi}\) is defined by the equation:
\begin{equation*}
(\,\widetilde d\,\Phi\,)_{(A,\lambda)}\,=\,\sigma_{(A,\lambda)} (X_{\Phi}\,,\,\cdot\,\,).
\end{equation*}
Now the partial derivative of \(\Phi\) at the point \((A,\lambda)\) to the direction \(a\in T_{A}{\cal A}\) is given by:
\[ \langle \,\frac{\delta \Phi}{\delta A}(A,\lambda)\,,\,a\,\rangle\,=\lim_{t\longrightarrow 0}\,\frac{1}{t}(\Phi(A+ta,\lambda)-\Phi(A,\lambda)).\]
Hence we have \(\frac{\delta \Phi}{\delta A}(A,\lambda)\,\in T^{\ast}_A{\cal A}\).      
Similarly the partial derivative of \(\Phi\) at the point \((A,\lambda)\) to the direction \(\alpha\in T^{\ast}_A{\cal A}\) is  given by 
\[ \langle \,\alpha\,,\,\frac{\delta \Phi}{\delta \lambda}(A,\lambda)\,\rangle\,=\lim_{t\longrightarrow 0}\,\frac{1}{t}(\Phi(A,\lambda+t\alpha)-\Phi(A,\lambda)).\]
Hence \(\frac{\delta \Phi}{\delta \lambda}(A,\lambda) \in T_A{\cal A}\).   
It holds that 
\begin{equation*}
(\,\widetilde d\,\Phi\,)_{(A,\lambda)}\left(\begin{array}{c}a\\ \alpha\end{array}\right)\,=\,\langle\,
\frac{\delta \Phi}{\delta A}(A,\lambda),\,a\,\rangle\,+\,\langle \alpha\,,\,\frac{\delta \Phi}{\delta \lambda}(A,\lambda)\,\rangle\,,
\end{equation*}
for \(\forall\,
\left(\begin{array}{c}a\\ \alpha\end{array}\right)\in T_{(A,\lambda)}T^{\ast}{\cal A}\).
 
So the Hamiltonian vector field of \(\Phi\) is given by
\begin{equation*}
X_{\Phi}\,=\,\left(\begin{array}{c}-\frac{\delta \Phi}{\delta \lambda}\\[0.2cm]
\frac{\delta \Phi}{\delta A}\end{array}\right),
\end{equation*}
\cite{MW}.
For example, if we take the Hamiltonian function 
\begin{equation*}
H(A,\lambda)\,=\, \frac12\int_M\,tr[F_A\wedge\ast F_A\,]+\,\frac12\int_M\,tr[\lambda\wedge \ast\lambda]\,,
\end{equation*}
then
\begin{equation*}
X_{H}\,=\,\left(\begin{array}{c} -\,\ast\lambda \\[0.2cm]
 d_A(\ast F_A) \end{array}\right).
\end{equation*}
It follows that the critical points of the Hamiltonian function \(H(A,0)= \frac12\int_Mtr[F_A\wedge\ast F_A\,]\)  are given by the Yang-Mills equation: \(d_AF_A=d_A^{\ast}F_A=0\).   
\\

The group of gauge transformations \({\cal G}(M)=\mathcal{E}^0_l(M,\,Lie\,G)\) acts on \(T_A{\cal A}\)  by the adjoint representation; 
\(a\longrightarrow Ad_{g^{-1}}\,a=g^{-1}ag\), and on  \(T^{\ast}_A{\cal A}\) by its dual \(\alpha\longrightarrow g\alpha g^{-1}\).     Hence the canonical 1-form and 2-form are \({\cal G}\)-invariant.   
The infinitesimal action of \(\xi\in Lie\,{\cal G}\) on the cotangent space \(T^{\ast}{\cal A}\)  
gives a vector field \(\xi_{T^{\ast}{\cal A}}\), called the {\it fundamental vector field on \(T^{\ast}{\cal A}\)}, that is defined at the point \((A,\lambda)\) by the equation:
\begin{equation*}
\xi_{T^{\ast}{\cal A}}(A,\lambda)=
\frac{d}{dt}\,\exp\,t\xi\cdot  \left(\begin{array}{c}A\\[0.2cm]
 \lambda\end{array}\right)=\left(\begin{array}{c}d_A\xi\\ [0.2cm]
\left[\xi,\lambda\right]\end{array}\right).
\end{equation*}

The moment map of the action of \(\,{\cal G}\) on the symplectic space \(\,(T^{\ast}{\cal A},\,\sigma)\) can be  described as follows.   For each \(\xi\in Lie\, {\cal G}\) we define the function   
\begin{equation}\label{momentmap}
J^{\xi}(A,\lambda)=\,\theta_{(A,\lambda)}\left(\xi_{T^{\ast}{\cal A}}\right)\,=\int_M\,tr\,(\,d_A\xi\wedge\,\lambda)\,.
\end{equation} 
Then the correspondence \(\xi\longrightarrow J^{\xi}(A,\lambda)\) is linear and defines a element of \(
J(A,\lambda)\in Lie\,{\cal G}^{\ast}\) and we have a map
\begin{eqnarray*}
& J\,: \,T^{\ast}{\cal A}\ni (A,\lambda)\,\longrightarrow J(A,\lambda)\in (Lie\,{\cal G})^{\ast}\,\simeq\,\mathcal{E}^{m-1}(M,\,Lie\,G),\\[0.2cm]
& J^{\xi}(A,\lambda)=\,<\,J(A,\lambda\,)\,,\,\xi\,>.
\end{eqnarray*}

(\ref{momentmap}) yields 
\begin{equation*}
\widetilde{d}\,J^{\xi}\,=\,\sigma(
\, \xi_{T^{\ast}{\cal A}}\,\,,\quad\cdot\quad)\,,\quad\mbox{ for }\,\forall \xi\in Lie\,{\cal G}.
\end{equation*}

Hence we have the following: 
\begin{thm}
\begin{enumerate}
\item
The action of the group of gauge transformations \({\cal G}(M)\) on the symplectic space \((\,T^{\ast}{\cal A}(M),\,\sigma\,)\) is a Hamiltonian action and the moment map is given by 
\begin{equation*}
  J^{\xi}(A,\lambda)=\,\int_M\,tr\,(\,d_A\xi\wedge\,\lambda)\,,\quad\,\forall \xi\in {\cal G}(M).
\end{equation*} 
\item
In the case when \(M\) has boundary, 
the action of the group of gauge transformations \({\cal G}_0(M)\) on the symplectic space \((\,T^{\ast}{\cal A}(M),\,\sigma\,)\) is a Hamiltonian action and the moment map is given by 
\begin{equation}
J_0(A,\lambda)=\,-\, d_A\,\lambda\,. \label{moment}
\end{equation} 
\end{enumerate}
\end{thm}
The second assertion follows from Stokes' theorem since any \(\xi\in Lie\,{\cal G}_0\) has the boundary value \(0\):
\[J_0^{\xi}(A,\lambda)=\,\int_M\,tr( d_A\xi\wedge\lambda)=\,-\,\int_M\,tr( \xi\wedge\,d_A\lambda)\,.\]

\subsection{Generating functions}

Suppose that we are given a section of \(T^{\ast}{\cal A}\):
\begin{equation*}
\phi\,:\,{\cal A}\,\longrightarrow\,T^{\ast}{\cal A}.
\end{equation*}
We denote it  \(\,\phi(A)=(A,\,f(A))\) with \(\,f(A)\in T_A^{\ast}{\cal A}\).     

Let \( {\phi}^{\ast}\theta\) be the pullback by \(\phi\) of the canonical 1-form \(\theta\) on \(T^{\ast}\mathcal{A}\).   We have the tautological relation:
\(\,{\phi}^{\ast}\theta=f\).
That is,
\begin{equation*}
({\phi}^{\ast}\theta\,)_A(a)=\,\langle \,f(A),\,a\rangle\,,\quad \forall a \in T_A{\cal A}\,.  
\end{equation*}

Let \( \phi^{\,\ast}\sigma\) be the pullback by \(\phi\) of the canonical 2-form \(\sigma\):
\begin{equation*}
(\phi^{\,\ast}\sigma)_A(a,b)=\sigma_{\phi(A)}(\phi_{\ast}a,\,\phi_{\ast}b)=
\sigma_{(A,f(A))}(\left(\begin{array}{c}a\\ (f_{\ast})_Aa\end{array}\right),\left(\begin{array}{c}b\\ (f_{\ast})_Ab\end{array}\right)\,).
\end{equation*}
\( \phi^{\,\ast}\sigma\)  is a closed 2-form on \({\cal A}\).     The relation \({\phi}^{\ast}\theta=f\) implies the relation
\begin{equation*}
\phi^{\,\ast}\sigma=\widetilde d\,f\,.
\end{equation*}   
\(f\) is called a (local) {\it generating function}.   
\(\phi^{\,\ast}\sigma\) is not necessarily non-degenerate.    \\

{\bf Example}[Atiyah-Bott, 1982],\cite{AB}.\\
Let \(\Sigma\) be a surface ( 2-dimensional manifold ).    \\
\[T_A{\cal A}(\Sigma) \simeq\,T_A^{\ast}{\cal A}(\Sigma)\simeq  \mathcal{E}^1_1(\Sigma, Lie G).\]
Define the section  \( \phi(A)=(A,f(A))\) of the cotangent bundle \(T^{\ast}{\cal A}\) by the generating function 
\[ f:\,{\cal A}\ni A \longrightarrow f(A)=A\in   \mathcal{E}^1_1(\Sigma, Lie G)=T_A^{\ast}{\cal A}(\Sigma) .\]  
Then  
\[({\phi}^{\ast}\theta)_Aa=\int_\Sigma\,tr(\,Aa\,)\,.\]
We have the 2-form;
\begin{eqnarray*}
(\phi^{\,\ast}\sigma)_A(a,b)&=&\,(\widetilde d\,( {\phi}^{\ast}\theta))_A(a,b)\,=\,
\langle\, (\partial {\phi}^{\ast}\theta(A))a,b\,\rangle - \langle \,(\partial {\phi}^{\ast}\theta(A))b, a\,\rangle \\[0.2cm]
&=&\,\int_\Sigma\,tr(ba)\,-\,\int_\Sigma\,tr(ab)=2\int_\Sigma\,tr(ba).\end{eqnarray*}

 \(({\cal A}(\Sigma),\,{\phi}^{\ast}\sigma\,)\) is a symplectic manifold, in fact \({\phi}^{\ast}\sigma\) is non-degenerate.\\

\section{Pre-symplectic structure on the space of connections on a four-manifold}

From this section we shall deal with the space of connections {\it on a four-manifold and on its three dimensional boundary}.     
Let \(X\) be an oriented Riemannian four-manifold with the boundary \(M=\partial X\) that may be empty.      

We consider the trivial principal bundle:  \(P=X\times SU(n)\),\(\, n\geq 3\).     We denote the group of \(L^2_3\)-gauge transformations on \(P\) by \(\mathcal{G}(X)\).    The subgroup of those transformations that are the identity on the boundary \(M\) is denoted by \(\mathcal{G}_0(X)\).   We denote the space of irreducible \(L^2_{2}\)-connections by \({\cal A}(X)\) which is abbreviated to \({\cal A}\) when there is no confusion.   The tangent space at \(A\in {\cal A}\) is
  \[T_A{\cal A}\,=\,\mathcal{E}^1_{2}(X,Lie\,G).\] 
  Hence the cotangent space at \(A\in \mathcal{A}\) is
  \[T^{\ast}_A{\mathcal{A}}\,=\,\mathcal{E}^3_{2}(X,Lie\,G).\]

   We define a section \( \vartheta\) of the cotangent bundle  \(T^{\ast}\mathcal{A}\,\stackrel{\pi}{\,\longrightarrow\,} {\cal A}\) by 
\begin{equation*}
\vartheta(A)\,=\left(\,A\,,\,q(\,AF_A+F_AA-\frac12 A^3)\right)\,.
\end{equation*}
The generating function \(q(AF_A+F_AA-\frac12 A^3)\) is called the \textit{Chern-Simons 3-form} on \(X\) valued in \(su(n)\), where \(q=\frac{1}{24\pi^3} \).\\

The tangential map of \(\,\vartheta:\,\mathcal{A}\longrightarrow T^{\ast}\mathcal{A}\)  at \(A\in \mathcal{A}\) becomes 
\begin{equation}\label{tgCS}
(\vartheta_{\ast})_Aa=\left(\begin{array}{c}a\\ [0.2cm]
q(aF_A+F_Aa+Ad_Aa+d_AaA-\frac12(aA^2+AaA+A^2a))\end{array}\right),
\end{equation}
for any \(a\in T_A{\cal A}\,\).

Let \(\theta\) be the canonical 1-form  (\ref{canonical1}) and let \(\sigma\) be the canonical 2-form (\ref{canonical2}) on \(T^{\ast}\mathcal{A}\).

\begin{lem}
Let \(\,\Theta=\vartheta^{\ast}\,\theta\) and \(\,\Omega=\vartheta^{\ast}\,\sigma\)  be the differential forms on \(\mathcal{A}\) that are the pullbacks by \(\vartheta\) of the canonical forms \(\theta\) and \(\sigma\)  .   Then we have 
\begin{equation*}
\Theta_A(\,a\,)=\frac{1}{24\pi^3}\int_X\,tr[\, (A F+F A-\frac12 A^3)\, a \,]\,,
\end{equation*}
\begin{equation}\label{symp2form}
\Omega_A(\,a,\,b\,)=\frac{1}{8\pi^3}\int_Xtr[(ab-ba)F ] -\frac{1}{24\pi^3}\int_{M}tr[(ab-ba)A]\,.
\end{equation}
for \( a,\,b\,\in T_A{\cal A}\).
\end{lem}
 
The first equation follows from the tautological relation of the canonical 1-form \(\theta\).   
For \(a,b \in T_A{\cal A}\),  we have
\begin{eqnarray*}
\Omega_A(\,a,\,b\,)&=&
(\widetilde d\,\Theta)_A(\,a,\,b\,)=\langle\, (\partial\,\Theta(A))a\,,\,b\,\rangle - \langle\, (\partial\,\Theta(A))b\,, \,a\,\rangle \\[0.2cm]
&=&  \frac{1}{24\pi^3}\int_X\,tr[\,2(ab-ba)F-(ab-ba)A^2\\[0.2cm]
&& \qquad\qquad - \,(bd_Aa+ d_Aab-d_Aba-ad_Ab)A\,] .
\end{eqnarray*}
Since
\[d\,tr[(ab-ba)A]=tr[(b\,d_Aa+ d_Aa\,b-d_Ab\,a-a\,d_Ab)A]+tr[(ab-ba)(F+A^2)],\]
we have 
\begin{equation*}
\Omega_A(\,a,\,b\,)=\frac{1}{8\pi^3}\,\int_X\,tr[\,(a b-b a)F\,] \,-\,\frac{1}{24\pi^3}\int_{M}\,tr[\,(a b-b a) A\,]\,, 
\end{equation*}
for  \(a,b\in T_A{\cal A}\).   
\hfill\qed

Thus we have the following:
\begin{thm}\cite{K}\label{presymp}~~~
Let \(P=X\times SU(n)\) be the trivial \(SU(n)-\)principal bundle over a four-manifold \(X\).   There exists a pre-symplectic structure on \({\cal A}(X)\) given by 
the 2-form (\ref{symp2form}).
\end{thm}

Since any \(g\in \mathcal{G}_0(X)\) is the identity transformation on boundary \(M\), the 2-form \(\Omega\,\) is \({\cal G}_0(X)\)-invariant, so by the Marsden-Weinstein reduction theorem \cite{MR} applied to the pre-symplectic space \(({\cal A}\,,\,\Omega )\)  we have the following:
\begin{cor}~~
There exists a pre-symplectic structure on the moduli space of connections \(\,{\cal B}(X)=\mathcal{A}(X)/\mathcal{G}_0(X)\).
\end{cor}

If \(X\) has no boundary and \(A\) is a flat connection then \(\Omega=0\), so we have the following
\begin{cor}~~~
Let \(X\)  be a compact 4-manifold without boundary then
\[L^{CS}=\{\,\vartheta(A)\,;\quad A\in {\cal A}^{\flat}(X)\} \,\] is a 
Lagrangian submanifold of \(T^{\ast}{\cal A}(X)\).
\end{cor}

In fact the tangential map \(\vartheta_{\ast}\) of \(\vartheta\), (\ref{tgCS}), is an isomorphism of \(T\mathcal{A}^{\flat}\), so \(L^ {CS}\) becomes a submanifold of \(T^{\ast}{\cal A}\).

\begin{rem}~~
\(\Theta\,\) is not \({\cal G}_0(X)\)-invariant.   In fact we gave in  \cite{K} the variational formula  of \(\Theta\) under the action of \(\,{\cal G}_0(X)\) that was important to construct an hermitian line bundle with connection over \({\cal B}(X)\). 
\end{rem}

 Since \(\vartheta\) is not \({\cal G}_0(X)\)-equivariant, the pullback \(\vartheta^{\ast}J_0\) of the moment map \(J_0:T^{\ast}{\cal A}\longrightarrow (Lie\,{\cal G}_0)^{\ast}\), (\ref{moment}),  does not give a moment map on \(\mathcal{A}(X)\) under the action of  \({\cal G}_0(X)\).     Instead we have the following theorem.
 
\begin{thm}\cite{K}~~
The action of  \({\cal G}_0(X)\) on \({\cal A}(X)\) is a Hamiltonian action with the corresponding moment map given by 
\begin{eqnarray*}
&&\Phi:\,{\cal A}(X)\longrightarrow\,( Lie\,{\cal G}_0)^{\ast}\,;\, A\longrightarrow\,F_A^2\,,\\[0.2cm]
&&\langle \Phi(A),\,\xi\rangle\,=\Phi^{\xi}(A)\,=\, \frac{1}{8\pi^3}\int_X\,tr(\,F_A^2\,\xi\,)\,,\quad \forall\xi\in Lie\,{\cal G}_0\,.
\end{eqnarray*}
\end{thm}

{\it Proof}

The equivariance of \(\Phi:{\cal A}\longrightarrow (Lie\, {\cal G}_0)^{\ast}\) with respect to the \({\cal G}_0\)-action on \({\cal A}\,\) and the coadjoint action on \((Lie\, {\cal G}_0)^{\ast}\,\) is evident.     We have 
\[(\partial \Phi^{\xi}(A))\,a\,=\,\frac{1}{8\pi^3}\int_X\,Tr[(d_Aa\wedge F_A+F_A\wedge d_Aa)\xi\,],\]
and
\begin{equation}\label{moment1}
\begin{split}
&\int_X\,Tr[(d_Aa\wedge F_A+F_A\wedge d_Aa)\xi]\,-\,\int_X\,Tr[(a\,d_A\xi-d_A\xi \,a)F_A\,]\\[0.2cm]
&=\int_Xd\,Tr[(aF_A+F_Aa)\xi]=\int_{M}Tr[(aF_A+F_Aa)\,\xi\,]\,=\,0.
\end{split}
\end{equation}
The last equality follows from the fact that \(\xi=0\) on \( M\).   \\
   On the other hand 
 \begin{eqnarray*}
dTr[(Aa-aA)\xi\,]&=&Tr[(F_Aa+aF_A-Ad_Aa-d_AaA+A^2a+aA^2)\xi\,]\\[0.2cm]
&&\qquad \qquad \qquad +\, Tr[ (ad_A\xi-d_A\xi a)A\,]
,\end{eqnarray*}
the first term of which vanishes on \( M\).
 Hence  
\begin{equation}\label{btgorbit}
\int_MTr[ (a\,d_A\xi-d_A\xi \,a)A]
=\int_{M}
dTr[(Aa-aA)\xi]=\,0.\end{equation}
(\ref{moment1}) and (\ref{btgorbit}) imply 
\[ (\partial \Phi^{\xi}(A))a\,=\,\Omega_A(a,d_A\xi).\]
\hfill\qed 

~~\\

We have 
\(\Phi^{-1}(0)=\{A\in {\cal A}(X);\, F_A^2=0\}\) .      
Since \(\Omega\) is \({\cal G}_0\)-invariant and    \(i_{d_A\xi}\Omega=0\) on \(\Phi^{-1}(0)\) for \(\xi\in Lie\,{\cal G}_0\,\), (\ref{btgorbit}), 
the 2-form \(\Omega\) descends to \(\Phi^{-1}(0)/{\cal G}_0\) and gives a closed 2-form on it.   This is the  pre-symplectic reduction if \(\Phi^{-1}(0)\) is a manifold and the action of \(\,{\cal G}_0\) on it is locally free.         The 
\({\cal G}_0\)-action is indeed free but \(0\) is not a regular value of the moment map \(\Phi\), so we can not expect that \(\Phi^{-1}(0)\) is a smooth manifold.   Instead it contains the subspace of flat connections  \({\cal A}^{\flat}(X)\), which is known to be a smooth manifold if \( M=\partial X \neq\emptyset\) .     Therefore we have the following:
\begin{thm}\label{bsymplectic}~
Suppose \( M=\partial X\neq\emptyset\).    Then the moduli space of flat connections 
\({\cal M}^{\flat}(X)\), (\ref{flatmoduli}),  is a smooth manifold endowed with a pre-symplectic structure \(\omega^{\flat}\) that is given by 
\begin{equation*}
\omega^{\flat}_{[A]}(a, b)=\,-\frac{1}{24\pi^3}\int_{ M}tr[(ab-ba)A]
\end{equation*}
for \([A]\in{\cal M}^{\flat}(X)\) and \(a,\,b \in T_{[A]}{\cal M}^{\flat}(X)\).   
\end{thm}
 Here \([A]\in {\cal M}^{\flat}(X)\) denotes the \({\cal G}_0\)-orbit of \(A\in{\cal A}^{\flat}(X)\), and for a tangent vector \(a\in T_{[A]}{\cal M}^{\flat}(X)\) we take the representative tangent vector to the slice; \(a\in H_A^{\flat}\), (\ref{flattangentsp}).   In fact 
\(\,\omega^{\flat}\) is well defined by virtue of (\ref{btgorbit}).

\section{The space of flat connections on a three-manifold}

In this section we study the space of connections on a 3-manifold \(M\) by looking at the space of connections on a 4-manifold \(X\) whose boundary is \(M\); \(\partial X=M\).     

\subsection{Chern-Simons function}

It is a well known fact that given a principal \(G-\)bundle \(P\) over a 3-manifold \(M\) there exists   
an oriented 4-manifold \(X\) with boundary \(\partial X=M\) and a \(G-\)bundle \({\bf P}\) over \(X\) that extends  \(P\).    And any connection \(A\) on \(P\) has an extension to a   connection \({\bf A}\) on  
\({\bf P}\).    These are essentially consequences of Tietz's extension theorem for continuous functions on a closed subset \cite{KN}.

As before \(\mathcal{A}(X)\) denotes the space of irreducible \(L^2_{2}\)-connections on the principal bundle \(\mathbf{P}\longrightarrow X\), while \(\mathcal{A}(M)\) denotes the space of irreducible \(L^2_{\frac32}\)-connections on the principal bundle \(P\longrightarrow M\). 

We denote by
\begin{equation*}
r:\,{\cal A}(X)\longrightarrow\,{\cal A}(M)
\end{equation*}
the restriction  of connections on \(X\) to the boundary \(M\):
\[r(\mathbf{A})=A\vert M, \quad \mathbf{A}\in{\cal A}(X).\]
The tangent map of \(r\) at \({\bf A}\in{\cal A}(X)\) is
\[\rho_{ {\bf A}}\,:\,T_{\bf A}{\cal A}(X)=\mathcal{E}^1_{2}(X,Lie\,G)\longrightarrow T_A{\cal A}(M)=\mathcal{E}^1_{\frac32}(M,Lie\,G), \]
where \(A=r({\bf A})\).

 We generally use bold letters for connections on a manifold that extend connections on its boundary,  but sometimes we use plain letters when no confusion is likely.

Let \(\mathcal{G}(X)\) denote the group of \(L^2_3\)-gauge transformations on \(X\) and let \(\mathcal{G}(M)\) denote the group of \(L^2_{\frac52}\)-gauge transformations on \(M\).   
 The boundary restriction
 \[r:\, \mathcal{G}(X)=\mathcal{E}^0_3(X,G)\,\longrightarrow\,\mathcal{G}(M)=\mathcal{E}^0_{\frac52}(X,G)\] 
 is continuous.
 
  \( {\cal G}(M)\)  is not connected and decomposes into countably many sectors labeled by the mapping degree
\begin{equation*}
\deg f=\,\frac{1}{24\pi^2}\int_M\,Tr\,(df\,f^{-1})^3\,.\end{equation*}

We have the following relation:
\begin{equation*}
\deg(g\,f)=\deg(f)+\deg(g). 
\end{equation*}
 The  group of \(L^2_{3}\)-gauge transformations on \(X\) that are the identity on the boundary  \(M\) is denoted by \({\cal G}_0(X)\).   It is the kernel of the 
 restriction map \(r:{\cal G}(X)\longrightarrow{\cal G}(M)\).    If \(X\) is simply connected then 
\(f\in {\cal G}(M)\)  is the restriction to \(M\) of a \({\bf f}\in {\cal G}(X)\) if and only if  \(\deg f=0\).  Thus we have the following exact sequence:
\begin{equation*}
1\longrightarrow {\cal G}_0(X)\longrightarrow {\cal G}(X)\stackrel{r}{\longrightarrow } {\cal G}^{\deg 0}(M) \longrightarrow 1.
\end{equation*}
Here we denote 
\begin{equation}\label{deg0gauge}
 \mathcal{G}^{\deg 0}(M)=\{g\in {\cal G}(M);\quad \deg g=0\,\}.
 \end{equation}
 
From now on \(X\) is assumed to be a simply connected four-manifold with boundary \(M=\partial X\).   
On a 3-manifold any principal bundle has a trivialization.   We choose a trivialization so that every  connection is identified with a Lie algebra-valued 1-form.   We  define the degree of a connection by using the {\it 3-dimensional Chern-Simons function} \cite{D2}:
\begin{equation}\label{degA}
\deg(A)\,=- \frac{1}{8\pi^2}\int_{M}\,Tr\,(AF_A-\frac13A^3), \quad A\in {\cal A}(M).
\end{equation}
It depends on the trivialization only up to an integer.   When \(A\) is of pure-gauge, that is, \(A=f^{-1}df\) for a 
\(f\in \mathcal{G}(M)\), we have \(\deg(A)=\deg f\).   
From Stokes' theorem,  we have the well known relation: 
\begin{equation*}
\int_X\,Tr[\,F_{{\bf A}}^2]=\int_M\,Tr[AF_A-\frac13A^3]\,,\quad r(\mathbf{A})=A\,. 
\end{equation*}

The 3-dimensional Chern-Simons function descends to a map:
  \[\deg :\,{\cal B}(M)={\cal A}(M)/\mathcal{G}^{\deg 0}(M)\,\longrightarrow\,{\rm R}/{\bf Z}.\]
The critical points of the 3-dimensional Chern-Simons function are 
the gauge equivalence classes of flat connections on \(M\).   

\begin{prop}\label{degree}
For \(A\in{\cal A}(M)\) and \(g\in{\cal G}(M)\), we have 
\begin{equation*}
\deg (g\cdot A)=\deg (A)+\deg g.  
\end{equation*}
\end{prop}

\subsection{ A twisted pre-symplectic structure on the space of flat connections}

We introduce the following two differential forms on \(M\):
\begin{eqnarray*}
\omega_A(a,b)&=&-q\int_M\,Tr[(ab-ba)A],\\
\kappa_A(a,b,c)&=&-3q\int_M\,Tr[(ab-ba)c],\quad a,\,b,\,c \in T_A{\cal A}\,.\end{eqnarray*}

We have
\begin{equation*}
\widetilde d\,\omega_A=\kappa_A.\end{equation*}
In fact, for \(a,b,c\in T_A{\cal A}\), we have 
\[\widetilde d\,\omega_A(a,b,c)=3\partial_A(\omega_A(a,b))(c) =-3q\int_M\,Tr[(ab-ba)c]
=\kappa_A(a,b,c).\]
We call  \(({\cal A}(M),\,\omega, \kappa\,)\) a {\it pre-symplectic manifold twisted by the 3-form} \(\kappa\).

\begin{rem}~~~
For \(G=SU(2)\), 
\(\kappa\) and \(\omega\) vanish identically, \cite{Ko} Lemma 1.3.   So in the following we consider 
mainly the case \(G=SU(n)\) with \(n\geq 3\).
\end{rem}

Let \({\cal A}^{\flat}={\cal A}^{\flat}(M)\) be the space of flat connections;
 \[{\cal A}^{\flat}(M)=\{A\in{\cal A}(M);\,F_A=0\}.\]    

The tangent space of  \( {\cal A}^{\flat}\) at  \(A\in {\cal A}^{\flat}\)  is given by
\begin{equation*}
T_A{\cal A}^{\flat}=\{a\in\Omega^1(M,\,Lie\,G);\quad d_Aa=0\}.
\end{equation*}
 We have the orthogonal decomposition 
\begin{equation}
T_A{\cal A}^{\flat}=\{d_A\xi;\,\xi\in {\cal G}(M)\,\}\oplus H_A^{\flat},\label{3decomposition}
\end{equation}
where
\begin{equation*}
H^{\flat}_A= \{a\in \Omega^1(M,\,Lie\,G);\,d_A^{\ast}a=d_Aa=0\,\}.  \label{harmonic}
\end{equation*}
In fact \(\{d_A\xi;\,\xi\in \mathcal{G}(M)\}\) is tangent to \({\cal A}^{\flat}(M)\) because 
\(d_A(d_A\xi)=\,[\,F_A,\,\xi\,]\,=0\) on \({\cal A}^{\flat}(M)\).

 For \(A\in \mathcal{A}^{\flat}(M)\) let 
 \(H^1_{A}(M,\,Lie\,G)\) be the 1st cohomology class of the complex \(\,(\,\mathcal{E}^k(M,\,Lie\,G)\,,\,d_A\,)\).      We know from Hodge theory that it is equal to the space of harmonic 1-forms, so that 
\begin{equation*}
H^1_{A}(M,\,Lie\,G)\,=\, H^{\flat}_A\,.\end{equation*}

\begin{prop}
\({\cal A}^{\flat}(M)\) is a \({\cal G}(M)-\)invariant manifold, and  
 the action of \({\cal G}(M)\) on \({\cal A}^{\flat}(M)\) is infinitesimally symplectic, that is, the Lie derivative \(L_{d_A\xi }\,\omega\) vanishes on \({\cal A}^{\flat}(M)\).   
 \end{prop}
 
Evidently the condition \(F_A=0\) is \(\mathcal{G}(M)\)-invariant.    
The  proposition follows from the following lemma:
\begin{lem} \label{Lieder}
Let \(i_{d_A\xi}\) and \(L_{d_A\xi}\) denote respectively the inner derivative and the Lie derivative by the fundamental vector field \(d_A\xi\).    We have 
\begin{equation*}
 i_{d_A\xi}\,\kappa=0,\qquad L_{d_A\xi }\,\omega=0.
\end{equation*}
on  \( {\cal A}^{\flat}(M)\),
\end{lem}

{\it Proof}

We have, for \(a,b\in T_A{\cal A}^{\flat}\), 
\[ i_{d_A\xi }\kappa_A(a,b)=-3q\int_M\,Tr[(ab-ba)d_A\xi]=-3q\int_M\,dTr[(ab-ba)\xi]=0,
\]
because \(d_Aa=d_Ab=0\).
Then \( i_{d_A\xi }\widetilde d\,\omega=i_{d_A\xi }\kappa =0\) and 
\begin{equation*}
\begin{split}
(L_{d_A\xi}\omega)_A&(a,b)=(\widetilde d\,i_{d_A\xi}\,\omega)_A(a,b)
=\partial_A\,(i_{d_A\xi}\,\omega_A(b))(a)-\partial_A\,(i_{d_A\xi}\,\omega_A (a))(b)\\[0.2cm]
&= -\frac{1}{24\pi^3}\int_{ M}\,Tr[( b\,d_A\xi-d_A\xi\, b)a]+
\frac{1}{24\pi^3}\int_{ M}\,Tr[( a\,d_A\xi-d_A\xi \,a)b]\\[0.2cm]
&= -\frac{1}{12\pi^3}\int_{M}\,Tr[(ab-ba)d_A\xi]
= -\frac{1}{12\pi^3}\int_{ M}\,d\,Tr[(ab-ba)\xi]\\[0.2cm]
&=0,
\end{split}
\end{equation*}
for  \(A\in{\cal A}^{\flat}\) and for 
\(a,b\in T_A{\cal A}^{\flat}\).   
\hfill\qed 

\begin{rem}~~~
The 1-form \(i_{d_A\xi}\,\omega\) is explicitly given by 
\begin{equation*}
(i_{d_A\xi}\,\omega)_A\,(a)=q_3\,\int_M Tr[(A^2\xi+\xi A^2)a],
\end{equation*}
for \(a\in T_A{\cal A}^{\flat}(M)\).
\end{rem}

\subsection{Pre-symplectic sectors of the space of flat connections on \(M\)}

Let  
\begin{equation*}
r:\,{\cal A}^{\flat}(X)\,\longrightarrow \, {\cal A}^{\flat}(M)
\end{equation*}
be the restriction map to the boundary for the space of flat connections :
\[r({\bf A})={\bf A}\vert_M\,,\quad {\bf A}\in {\cal A}^{\flat}(X).\]
The tangent map of \(r\) at \({\bf A}\in {\cal A}^{\flat}(X)\) becomes
\[\rho_{\bf A}:\,\{{\bf a}\in \mathcal{E}^1(X, Lie\,G);\, d_{{\bf A}}{\bf a}=0\} \longrightarrow 
\{{ a}\in \mathcal{E}^1(M, Lie\,G);\, d_{A}{ a}=0\}.\]

We consider the Chern-Simons function \(\deg \) on \(\mathcal{A}^{\flat}(M)\), (\ref{degA}).     
The exterior derivative of \(\deg\) at \(A\in\mathcal{A}^{\flat}(M)\) becomes;
\begin{equation*}
\,( \widetilde d\,\deg)_A\,a\,= \frac{1}{8\pi^2}\, \int_{M}\,Tr(\,A^2a)
=\frac{1}{8\pi^2}\, \int_{M}\,dTr(Aa)=0.
\end{equation*} 
Hence  \(\deg \) is constant on every connected component of \({\cal A}^{\flat}(M)\).

 \begin{defn}
 For each \(k\in {\bf Z}\) we define 
 \begin{equation*}
{\cal A}^{\flat}_k(M)=\left\{A\in{\cal A}^{\flat}(M);\quad \deg(A)=k\,\right\}.
\end{equation*}
We call \({\cal A}^{\flat}_k(M)\) the \(k\)-sector of the space of flat connections.
\end{defn}
Proposition \ref{degree} yields the invariance of  \(\mathcal{A}_k^{\flat}(M)\)  under the action of \(\mathcal{G}^{\deg 0}(M)\), (\ref{deg0gauge}).

Any connection on \(M\) is the restriction of a connection on \(X\).    But it is not clear whether we can take a flat connection on \(X\) for a given flat connection on \(M\).    In the following we shall define the range of \(r:\,{\cal A}^{\flat}(X)\,\longrightarrow \, {\cal A}^{\flat}(M)\).   

Let \(A\in{\cal A}^{\flat}(M)\).    Let  \(\widetilde X\) be the universal covering of \(X\) and 
 \(\widetilde M\) be the  subset of \( \widetilde X\) that lies over \(M\).   Let \(f_A\) be the parallel transformation by \(A\) along the paths starting  from \(m_0\in M\).   It defines a  smooth map on the covering space  \(\widetilde M\); \(f=f_A\in Map (\widetilde M, G) \), such that \(f^{-1}\,df=A\).   Then the   
degree of \(f\)  is equal to
 \begin{equation*}
\deg f=\frac{1}{24\pi^2}\int_{M}\,Tr\,A^3=\,\deg(A).
\end{equation*}
If \(\deg A\) vanishes, 
then there is a \({\bf f}\in{\cal G}(\widetilde X)\) that extends \(f\).   Therefore 
\({\bf A}={\bf f}^{-1}d{\bf f}\in {\cal A}^{\flat}(X)\) gives a flat extension of \(A\) to \(X\);  \(r({\bf A})=A\).      We have seen that a flat connection on \(M\) is the boundary restriction of a flat connection on \(X\) if \(\deg A=0\).

\begin{prop}\label{flatexrension}~~~
For any 4-manifold \(X\) with boundary \(M\) we have the following properties:
\begin{enumerate}
\item
The image of \(r\) is precisely  \({\cal A}^{\flat}_0(M)\) . 
 \item
\({\cal A}^{\flat}_k(M)\) is invariant under the action of \(\mathcal{G}^{\deg 0}(M)\) for any \(k\in \mathbf{Z}\) .
 \item
  \(\,d_A(\,Lie\, {\cal G}(M))\subset T_A{\cal A}^{\flat}_0(M)\).
   \item
 The action of the group of gauge transformations \({\cal G}(M)\) on \({\cal A}^{\flat}_0(M)\) is infinitesimally symplectic.
\end{enumerate}
\end{prop}

{\it Proof}

It follows from the above discussion that any \(A\in{\cal A}^{\flat}_0(M)\) is the boundary restriction of a \({\bf A}\in{\cal A}^{\flat}(X)\).   Conversely let \(A=r({\bf A})\) for a 
 \({\bf A}\in{\cal A}^{\flat}(X)\).   Then 
 \[\int_MTr\,A^3=\int_X\,Tr\,{\bf A}^4=0,\]
 so \(A\in {\cal A}^{\flat}_0(M)\).   Thus, for any 4-manifold \(X\) with boundary \(M\) the image of \(r\) is precisely \({\cal A}^{\flat}_0(M)\).  The 2nd assertion follows from (\ref{degA}).    The 3rd and 4th assertions are restatement of the facts
 \[d_A\xi\in T_A{\cal A}^{\flat}(M), \quad L_{d_A\xi}\,\omega=0.\] 
\hfill\qed 

\vspace{0.3cm}

The orthogonal complement \(H^{\flat}_A(M)\) of \(d_A(\,Lie\,{\cal G}(M)) \) in \(T_A{\cal A}^{\flat}(M)\), (\ref{3decomposition}),  is the space of harmonic 1-forms on \(M\), so is identified with the cohomology class  \(H^1_A(M, Lie\, G)\).   This is non-zero if and only if the connection can be deformed infinitesimally within \({\cal M}^{\flat}(M)\).   

\begin{lem}~~~
Let \(X\) be a 4-manifold with boundary \(M\).    Then \(r \) is a submersion.
\end{lem}
Take \(A\in{\cal A}^{\flat}_0(M)\) and \( {\bf A}\in{\cal A}^{\flat}(X)\) such that 
\(r({\bf A})=A\).    Let \(a\in T_A{\cal A}^{\flat}(M)\).    
From (\ref {3decomposition}), \(a \) is decomposed into 
 \[a=d_A\xi+b,\]
 by  \(\xi\in \mathcal{E}^0(M,Lie\,G) \) and  \(b\in H^{\flat}_A(M)\).   Let \(\eta\in \mathcal{E}^0(X,Lie\,G)\) be an extension of \(\xi\) to \(X\), then 
 \[\rho_{\bf A}(d_{\bf A}\eta)=d_A\xi.\]
 On the other hand the spaces of \(\Delta_A\)-harmonic 1-forms \(H^{\flat}_A(M)\) and \(H^{\flat}_{{\bf A}}(X)\) are isomorphic to the cohomology group \(H^1_{ A}(M,Lie\,G)\) and \(H^1_{ {\bf A}}(X,Lie\,G)\) respectively.    Since the  cohomology groups with compact support in \(X\); \(H^k_{{\bf A},\,c}(X,Lie\,G)\),  vanish for \(k=1,2 \), we have
 \[H^1_{{\bf A}}(X,\,Lie\,G)\simeq H^1_{A}(M,\,Lie\,G).\]
Hence there is a \({\bf b}\in H^1_{{\bf A}}(X,\,Lie\,G)=H^{\flat}_{{\bf A}}(X)\) such that 
\[b=\rho_{\bf A}{\bf b}+d_A\alpha=\rho_{\bf A}({\bf b}+d_{\bf A}\beta),\]
with \(\beta\in \mathcal{E}^0(X,Lie\,G)\) and  \(\alpha=r_{X}(\beta)\in \mathcal{E}^0(M,Lie\,G)\).
\(({\bf b}+d_{\bf A}\beta)\) being in \(T_{\bf A}{\cal A}^{\flat}(X)\) the lemma is proved; 
\[\rho_{\bf A}(d_{{\bf A}}\eta+{\bf b}+d_{\bf A}\beta)=d_A\xi+b=a.\]
\hfill\qed

 \begin{thm}
\(\,({\cal A}^{\flat}_0(M) , \omega)\) is a pre-symplectic manifold.
\end{thm}
 
 {\it Proof}
 
  We must show  
 \[\widetilde d\,\omega_A=\kappa_A=0,\]
 for any \(A\in{\cal A}^{\flat}_0(M)\).     Let \(X\) be a 4-manifold with boundary 
 \(M\).    There is a \(G-\)bundle  \({\bf P}\) over \(X\) and a flat connection \({\bf A}\in \mathcal{A}^{\flat}(X)\) such that  \(A=r{\bf A}\),  Proposition \ref{flatexrension}-1.     
Let \(a,b,c\in T_A{\cal A}^{\flat}(M)\).   \(\rho_{\bf A}\) being surjective, there are \({\bf a},{\bf b},{\bf c}\in T_{\bf A}{\cal A}^{\flat}(X)\) that extend \(a, b, c\) respectively.   Then we have 
\begin{eqnarray*}
\kappa_A(a,b,c)&=&-q\int_M\,Tr[(ab-ba)c]\nonumber \\[0.2cm]
&=&-
q\int_X\,Tr[\,(d_{{\bf A}}{\bf a}\,{\bf b}-{\bf a}\,d_{{\bf A}}{\bf b}\,-\,d_{{\bf A}}{\bf b}\,{\bf a}\,+\,{\bf b}\,d_{{\bf A}}{\bf a})\,{\bf c}\,
+\,({\bf a}{\bf b}-{\bf b}{\bf a})\,d_{{\bf A}}{\bf c}\,]\nonumber \\[0.2cm]
&=&0,\label{eq:symp}
\end{eqnarray*}
because of \(d_{{\bf A}}{\bf a}=0\), etc.
\hfill\qed 

 \vspace{0.3cm}
 
 \subsection{Relation to the moduli space  of flat connections over the four-manifold}
 
 Let \({\cal M}^{\flat}(X)\) be ( as in 1.3 ) the moduli space  of flat connections over \(X\).   By  Theorem \ref{bsymplectic}  \({\cal M}^{\flat}(X)\) is endowed with the pre-symplectic structure  
\begin{equation}
 \omega^{\flat}_{[{\bf A}]}({\bf a},\,{\bf b})
 =-q\int_M\,Tr[(ab-ba)A]\,,\label{defomeg}
\end{equation}
for \({\bf A}\in{\cal A}^{\flat}(X)\) and \({\bf a},\,{\bf b}\in  T_{\bf A}{\cal A}^{\flat}(X)\), where 
\(A=r({\bf A})\) and \(a=\rho_ {\bf A}({\bf a})\),  \(b=\rho_ {\bf A}({\bf b})\).      The right hand side coincides with the pre-symplectic form \(\omega_A(a,b)\) on \({\cal A}^{\flat}_0(M)\).   

 Evidently we have \(r(g\cdot {\bf A})=r({\bf A})\) for \(g\in{\cal G}_0\) .    Hence \(r\) induces the map  
\begin{equation*}
\overline r\,:\,{\cal M}^{\flat}(X)\,\longrightarrow \, {\cal A}^{\flat}(M).
\end{equation*}

\begin{prop}
\(\overline r\) gives a diffeomorphism of \({\cal M}^{\flat}(X)\) to \({\cal A}^{\flat}_0(M)\).
\end{prop}

{\it Proof}\\
We have already seen that \(r: {\cal A}^{\flat}(X)\longrightarrow {\cal A}^{\flat}(M)\) is a surjective submersion.   Hence it is enough to prove that \(\overline r\) is injective immersion.    In fact, 
let  \({\bf A}_1,\,{\bf A}_2\in{\cal A}^{\flat}(X)\) be such that \(r({\bf A}_1)=r({\bf A}_2)\), and let 
  \(f_{{\bf A}_i}\), \(i=1,2\), be the parallel transformations by \({\bf A}_i\), \(i=1,2\),  respectively,  along paths starting from \(m_0\in M\).   It defines a  smooth map on the universal covering space \(\widetilde X\stackrel{\pi}{\longrightarrow} X \); \(f_i =f_{{\bf A}_i}\in Map (\widetilde X, G) \), such that \(f_i^{-1}\,df_i={\bf A}_i\).   
Since \(r({\bf A}_1)=r({\bf A}_2)\) these parallel transformations coincide along the  paths that are contained in \(M\), that is, \(f_1\) and \(f_2\) coincide on  the covering space \(\widetilde M=\pi^{-1}(M)\) of \(M\).
Let  
 \(\tilde g\in Map(\widetilde X, G)\) be such that \(f_2=\tilde g\cdot  f_1\).    Then  \(\tilde g\) descends to a \(g\in{\cal G}_0(X)\) such that \({\bf A}_2=g\cdot {\bf A}_1\).   Therefore  \( \overline r\) is injective. \\
 The restriction of \(d_{{\bf A}}\,Lie\,{\cal G}_0(X)\) on the boundary \(M\) is obviously \(0\).   
The orthogonal complement 
of  \(d_{{\bf A}}\,Lie\,{\cal G}_0(X)\) in \(T_{{\bf A}}{\cal A}^{\flat}(X)\) consists of those 
\({\bf a}\in \mathcal{E}^1(X, Lie\, G)\) that satisfies \(d_{{\bf A}}{\bf a}=d_{{\bf A}}^{\ast}{\bf a}=0\).  We know that a harmonic 1-form vanishing on the boundary must vanish.  Therefore \({\bf a}=0\) if \({\bf a}\vert M=0\).    We have proved  
\[\ker \rho_{{\bf  A}}=d_{{\bf A}}\,Lie\,{\cal G}_0(X).\]
Thus \( \overline r\) is an injective immersion.
\hfill\qed 

From (\ref{defomeg}) and the above theorem we have the following
\begin{thm}\label{symplectomo}
\[ \overline r\,: {\cal M}^{\flat}(X)\longrightarrow  {\cal A}^{\flat}_0(M)\]
gives an isomorphism of  pre-symplectic manifolds;
\begin{equation*}
\left( {\cal M}^{\flat}(X),\, \omega^{\flat} \right)\,\simeq \,\left({\cal A}_0^{\flat}(M),\omega\right).
\end{equation*}
\end{thm}
~~~

Finally we shall investigate the relation between the moduli spaces \(\mathcal{N}^{\flat}(X)\) and  \(\mathcal{M}^{\flat}(X)\).   

The group of gauge transformations \({\cal G}(X)\) acts on  \({\cal A}(X)\) and its restriction to the space \({\cal A}^{\flat}(X)\) becomes infinitesimally symplectic.    This is seen by exactly the same calculation as  in Lemmas \ref{Lieder} where it is proved that the action of \(\mathcal{G}(M)\) on \({\cal A}^{\flat}(M)\) is infinitesimally symplectic.   
Since 
\[{\cal N}^{\flat}(X)={\cal A}^{\flat}(X)/{\cal G}(X)\simeq{\cal M}^{\flat}(X)/\,\mathcal{G}^{\deg 0}(M),\quad
(\ref{flatmodulin}), \]
 we have the presymplectic reduction \(({\cal N}^{\flat}(X), \,\omega^{\flat})\) .   
Theorem \ref{symplectomo} implies the following equivalence of the moduli spaces of flat connections on \(X\) and \(M\).
\begin{prop}
\begin{equation*}
{\cal N}^{\flat}(X)\simeq\, {\cal A}^{\flat}_0(M)/\,\mathcal{G}^{\deg 0}(M).
\end{equation*}
\end{prop}

{\bf Acknowledgements}

I would like to express my thanks to Professor Martin Guest of Waseda University for his valuable comments on this work.

  \end{document}